\documentclass[12pt,a4wide]{article}
\usepackage{amsmath}
\usepackage[all]{xy}
\usepackage[mathscr]{euscript}
\usepackage{amssymb}
\usepackage{amscd}
\usepackage[dvips]{epsfig}
\usepackage{amsmath, amscd, amssymb}

\newtheorem{lemma}{Lemma}[section]
\newtheorem{theorem}[lemma]{Theorem}
\newtheorem{prop}[lemma]{Proposition}

\newtheorem{cor}[lemma]{Corollary}
\newtheorem{remark}[lemma]{Remark}
\newtheorem{Def}[lemma]{Definition}
\newtheorem{Ex}[lemma]{Example}
\newtheorem{conj}[lemma]{Conjecture}

\newcommand{\bma}{\mathbb{A}}
\newcommand{\x}{\mathbf{x}}
\newcommand{\vv}{\mathbf{v}}
\newcommand{\odm}{\bma \mathbf{x}}
\newcommand{\Ev}{{\rm{Spec} \,}}
\newcommand{\spn}{{\, \rm{span} \,}}
\newcommand{\rank}{{\, \rm{rank} \,}}

\newcommand{\diag}{{\, \rm{diag} \,}}

\newcommand{\epf}{\hfill\hbox{\rule{3pt}{6pt}}\\}

\begin{document}

\title{Characterizing completely regular codes from an algebraic viewpoint}

\author{J. H. Koolen\footnote{Pohang Mathematics Institute, POSTECH, South Korea}\hspace{1mm}$\hspace{1mm}^{,}$\hspace{1mm}\footnotemark[2]
\and W. S. Lee\footnote{Department of Mathematics, POSTECH, South Korea}
 \and W. J. Martin\addtocounter{footnote}{0}\footnote{
Department of Mathematical Sciences
and
Department of Computer Science,
Worcester Polytechnic Institute,
Worcester, Massachusetts, USA.}
}

\date{\today}

\maketitle

\begin{abstract}
Completely regular codes are rich substructures in distance-regular graphs and have been studied
extensively over the last two decades. The class includes highly structured and beautiful examples 
such as perfect and uniformly packed codes while the rich properties of these codes allow for both
combinatorial and algebraic analysis. In fact, these codes are fundamental to the study of distance-regular
graphs themselves.

In a companion paper, we study products of completely regular 
codes and codes whose parameters form arithmetic progressions. 
This family of completely regular codes, while quite special in 
one sense, contains some very important examples and exhibits 
some of the nicest features of the larger class. Here, we 
approach these features from an algebraic viewpoint, exploring 
$Q$-polynomial properties of completely regular codes.

We first summarize the basic structure of the outer distribution 
module of a completely regular code. Then, employing a simple 
lemma concerning eigenvectors in association schemes, we propose 
to study the tightest case, where the indices of the eigenspace 
that appear in the outer distribution module are equally spaced.
In addition to the arithmetic codes of the companion paper, this 
highly structured class includes other beautiful examples and we 
propose the classification of $Q$-polynomial completely regular 
codes in the Hamming graphs. A key result is Theorem 
\ref{TQiffLeonard} which finds that the $Q$-polynomial condition 
is equivalent to the presence of a certain Leonard pair. This 
connection has impact in two directions. First, the Leonard pairs 
are classified and we gain quite a bit of information about the 
algebraic structure of any code in our class. But also this gives 
a new setting for the study of Leonard pairs, one closely related 
to the classical one where a Leonard pair arises from each 
thin/dual-thin irreducible module of a Terwilliger algebra of some 
$P$- and $Q$-polynomial association scheme, yet not previously 
studied. It is particularly interesting that the Leonard pair 
associated to some code $C$  may belong to one family in the Askey 
scheme while the distance-regular graph in which the code is found 
may belong to another.
\end{abstract}

\section{Introduction}
\label{Sec:intro}
The study of digital error-correcting codes includes as an important
and intriguing sub-topic
the analysis and classification of highly regular codes. These include
the perfect codes as well as several phenomenal families such as the
Kerdock codes, the Delsarte-Goethals codes, and the Reed-Muller codes.
One motivation for this branch of coding theory has always been a well-studied
but 
mysterious connection to finite groups. Optimal codes tend to have a
great deal of symmetry (as is often true in optimization problems which
themselves are defined in a symmetric way), and several finite simple
groups -- namely the Mathieu groups -- play an important role in the classification
of perfect codes.

But the class of completely regular codes, which properly contains both
the class of perfect codes and the class of uniformly packed codes but also, for
example, the Preparata and Kasami codes,  has
not received a great deal of attention in recent years. Our view is that
these codes deserve further study, not only because of their
connection to highly symmetric codes and codes with large minimum distance,
but also because of a key role that completely regular codes play in the
study of distance-regular graphs. A theorem of Brouwer, et al.
\cite[p353]{B} states that every distance-regular graph on a prime power
number of vertices admitting an elementary abelian
group of automorphisms which acts transitively on its vertices
is a coset graph of some additive
completely regular code in some Hamming graph (with some conference graphs
as exceptions). This gives another reason why a careful study of completely regular codes
in Hamming graphs (and, more generally, in distance-regular graphs) is central
to the study of association schemes.

In a companion paper, we study products of completely regular 
codes and codes whose parameters form arithmetic progressions. 
This family of completely regular codes, while quite special in 
one sense, contains some very important examples and exhibits 
some of the nicest features of the larger class. Here, we 
approach these features from an algebraic viewpoint, exploring 
$Q$-polynomial properties of completely regular codes.

We first summarize the basic structure of the outer distribution 
module of a completely regular code. Then, employing a simple 
lemma concerning eigenvectors in association schemes, we propose 
to study the tightest case, where the indices of the eigenspace 
that appear in the outer distribution module are equally spaced.
In addition to the arithmetic codes of the companion paper, this 
highly structured class includes other beautiful examples and we 
propose the classification of $Q$-polynomial completely regular 
codes in the Hamming graphs. A key result is Theorem 
\ref{TQiffLeonard} which finds that the $Q$-polynomial condition 
is equivalent to the presence of a certain Leonard pair. This 
connection has impact in two directions. First, the Leonard pairs 
are classified and we gain quite a bit of information about the 
algebraic structure of any code in our class. But also this gives 
a new setting for the study of Leonard pairs, one closely related 
to the classical one where a Leonard pair arises from each 
thin/dual-thin irreducible module of a Terwilliger algebra of some 
$P$- and $Q$-polynomial association scheme, yet not previously 
studied. It is particularly interesting that the Leonard pair 
associated to some code $C$  may belong to one family in the Askey 
scheme while the distance-regular graph in which the code is found 
may belong to another.

\section{Preliminaries and definitions}
\label{Sec:prelim}

\subsection{Distance-regular graphs}
Suppose that $\Gamma$ is a finite, undirected, connected graph
with vertex set $V\Gamma$. For vertices $x$ and $y$ in $V\Gamma$,
let $d(x,y)$ denote the distance between $x$ and $y$,
i.e., the length of a shortest path connecting $x$ and $y$ in $\Gamma$.
Let $D$ denote the diameter of $\Gamma$; i.e., the maximal distance
between any two vertices in $V\Gamma$.
For $0\le i\le D$ and $x\in V\Gamma$, let
${\Gamma_{i}}(x):=\{y\in V \Gamma\mid d(x,y)=i\}$
and put
$\Gamma_{-1}(x):=\emptyset$, $\Gamma_{D+1}(x):=\emptyset$. The graph
$\Gamma$ is called {\em distance-regular} whenever it is regular
of valency $k$, and there are integers $b_{i},c_{i} \ (0\leq i \leq D)$
so that for
any two vertices $x$ and $y$ in $V\Gamma$ at distance $i$, there
are precisely $c_{i}$ neighbors of $y$ in $\Gamma_{i-1}(x)$ and
$b_{i}$ neighbors of $y$ in $\Gamma_{i+1}(x).$  It follows that
there are exactly $a_i=k-b_i-c_i$ neighbors of $y$ in $\Gamma_i(x)$.
The numbers $c_{i}$, $b_{i}$ and $a_{i}$ are called the
{\em intersection numbers} of $\Gamma$ and we observe that
$c_{0}=0$, $b_{D}=0$, $a_{0}=0$, $c_{1}=1$ and $b_{0}=k$. The array
$\iota(\Gamma):=\{b_{0},b_{1},\ldots,b_{D-1};c_{1},c_{2},\ldots,c_{D}\}$
is called the {\em intersection array} of $\Gamma$. Set the tridiagonal matrix
$$L(\Gamma):=\begin{pmatrix}
a_{0} & b_{0}  &    &         &  \\
c_{1} & a_{1} & b_{1} &  &   \\
   & c_{2} & a_{2} & b_{2}  &   \\
   &            &  \ddots & \ddots  &  \ddots  \\
   &           & & c_{D}& a_{D}\end{pmatrix}.$$

From now on, assume $\Gamma$ is a
distance-regular graph of valency $k\ge 2$ and diameter $D\ge 2$.
Define $A_{i}$ to be
the square matrix of size $|V\Gamma|$  whose rows and columns are
indexed by $V\Gamma$ with entries
$${{(A_{i})}_{xy}= \begin{cases} 1 \hspace{4mm}\text{if} \hspace{2mm} d(x,y)=i\\
0 \hspace{4mm}\text{otherwise}\end{cases}}\hspace{3mm}(0\leq i\leq
D, \ x,y \in V\Gamma).$$

We refer to $A_{i}$ as the {\em $i^{\rm th}$ distance matrix} of $\Gamma$.
We abbreviate $A:=A_{1}$ and call this the {\em adjacency matrix} of $\Gamma$.
Since $\Gamma$ is distance-regular, we have for $2\le i \le D$
$$ A A_{i-1} = b_{i-2} A_{i-2} + a_{i-1} A_{i-1} + c_i A_i $$
so that $A_{i}=p_{i}(A)$ for some polynomial $p_{i}(t)$ of degree $i$.
Let $\bma$ be the {\em Bose-Mesner algebra}, the matrix algebra over
$\mathbb{C}$ generated by $A$.
Then $\dim \bma =D+1$ and
$\{A_{i}\hspace{1mm}|\hspace{1mm}0 \le i \le D\}$
is a basis for $\bma$. As $\bma$ is semi-simple and commutative,
$\bma$ has also a basis of pairwise orthogonal idempotents
$\left\{ E_{0}=\frac{1}{|V\Gamma|}J,E_{1},\ldots,E_{D} \right\}$.
We call these matrices the {\em primitive idempotents} of $\Gamma$.
As $\bma$ is closed under the entry-wise (or Hadamard) product  $\circ$,
there exist real numbers $q_{ij}^\ell$, called the {\em Krein parameters}, such that
\begin{equation}
\label{Ekrein}
E_i\circ E_j =\frac{1}{|V\Gamma|}\sum_{\ell=0}^{D} q_{ij}^\ell E_\ell \ \
( 0\leq i,j \leq D)
\end{equation}

We say the distance-regular graph $\Gamma$ is {\em $Q$-polynomial} with respect to
a given ordering
$E_0,E_1,\ldots, E_D$ of its primitive idempotents provided its Krein
parameters satisfy
\begin{itemize}
\item $q_{ij}^\ell =0$ unless $|j-i| \le \ell \le i+j$;
\item $q_{ij}^\ell \neq 0$ whenever $\ell=|j-i|$ or $\ell= i+j \le D$.
\end{itemize}

By an {\em eigenvalue} of $\Gamma$, we mean an eigenvalue of $A=A_1$. Since
$\Gamma$ has diameter $D$, it has at least $D+1$ eigenvalues; but since $\Gamma$
is distance-regular, it has exactly $D+1$ eigenvalues\footnote{See, for example
Lemma 11.4.1 in \cite{Godsil}.}, and they are exactly the eigenvalues of $L(\Gamma)$.

We denote these eigenvalues by
$\theta_0,\ldots,\theta_D$ and, aside from the convention that $\theta_0=k$,
the valency of $\Gamma$, we make no further assumptions at this point
about the eigenvalues except that they are distinct. We note that, with an appropriate
ordering of the eigenvalues, the
$i^{\rm th}$ primitive idempotent $E_i$ is precisely the matrix representing
orthogonal projection onto $V_i$, the eigenspace of $A$ associated to $\theta_i$.
In fact, when $\theta=\theta_i$ for some $i$, we will sometimes write $E(\theta)$
in place of $E_i$ when it is convenient to omit the subscript.

The following fundamental result will be very useful in this paper; it
is originally due to Cameron, Goethals, and Seidel \cite{C}.
\begin{theorem}[{\cite[Theorem~5.1]{C}}]
\label{Tcgs}
If ${\bf u}\in V_{i}$ and ${\bf v}\in V_{j}$ and $q_{ij}^{\ell}=0,$ then
${\bf u}\circ {\bf v}$ is orthogonal to $V_{\ell}$ where
${\bf u}\circ {\bf v}$ denotes the entry-wise product of vectors ${\bf u}$ and
${\bf v}$.\epf
\end{theorem}
An elementary proof of this fact can be found in \cite{martin-sym}.

For each eigenvalue $\theta$ of $\Gamma$ and for each $x\in V\Gamma$, there is
a unique normalized eigenvector in $V_i$ which is constant over each vertex
subset $\Gamma_i(x)$. The entries of this eigenvector, which we shall denote
by $u_i(\theta)$ ($0\le i\le D$, $\theta$ an eigenvalue of $\Gamma$) are
determined entirely by the intersection array, independent of the choice of $x$.

Suppose $\Gamma$ has intersection array
$\iota(\Gamma):=\{b_{0},b_{1},\ldots,b_{D-1};c_{1},c_{2},\ldots,c_{D}\}$.
and let $\theta$ be an eigenvalue of $\Gamma$. The
corresponding {\em standard right eigenvector }
$[u_{0}(\theta)=1,u_{1}(\theta),\ldots,u_{D}(\theta)]^\top$
of $\Gamma$ with respect to $\theta$ is defined by
the following initial conditions and recurrence relation:
\begin{equation}
\begin{split}
u_{0}(\theta)=1, \ u_{1}(\theta)=\theta/k,\hspace{30mm}\\
c_{i}u_{i-1}(\theta)+a_{i}u_{i}(\theta)+b_{i}u_{i+1}(\theta)=\theta u_{i}(\theta) \  \ (0\leq i \leq D),\\
\text{where}\hspace{2mm}u_{-1}=u_{D+1}=0.\hspace{30mm}
\end{split}
\end{equation}
One easily checks that the vector $\mathbf{u}$ of length $|V\Gamma|$
satisfying $\mathbf{u}_y = u_i(\theta)$
whenever $y \in \Gamma_i(x)$ satisfies $A\mathbf{u} = \theta \mathbf{u}$.  Let
$x\in V\Gamma$ and let $\mathbf{e}_x$ denote the elementary basis vector in $V$
corresponding to $x$.  For $\theta=\theta_j$ ($0\le j\le D$) we easily see that
$$ E_j \mathbf{e}_x = \frac{ m_j }{ |V\Gamma| } \mathbf{u}  $$
where $m_j := \rank E_j$.
It follows from this and (\ref{Ekrein})  that, for $0\le h,i,j \le D$,
\begin{equation}
\label{Ekreinu}
m_i m_j \, u_h (\theta_i) u_h( \theta_j ) = \sum_{\ell =0}^D q_{ij}^\ell
m_\ell  u_h(\theta_\ell).
\end{equation}
So we can detect whether or not $\Gamma$ is $Q$-polynomial just by looking
at its standard right eigenvectors.

\subsection{Codes in distance-regular graphs}

Let $\Gamma$ be a distance-regular graph with distinct eigenvalues
$\theta_{0}=k,\theta_{1},\ldots,\theta_{D}$.
By a {\em code} in $\Gamma$, we simply mean any
nonempty subset $C$ of $V \Gamma$.
We call $C$ {\em trivial} if $|C|\leq 1$ or $C=V\Gamma$ and {\em non-trivial} otherwise.
For $|C|>1$, the {\em minimum distance} of $C$, $\delta(C)$, is defined as
$$\delta(C):= \min \{\hspace{1mm}d(x,y) \ | \ x,y\in C,x\neq
y\hspace{1mm}\}$$ and
for any $x \in V \Gamma $ the distance $d(x,C)$ from $x$ to $C$ is defined as
$$d(x,C):= \min \{\hspace{1mm}d(x,y) \ | \ y\in C\hspace{1mm}\}.$$ The number
$$\rho(C):= \max \{\hspace{1mm}d(x,C) \ | \ x\in V
\Gamma\hspace{1mm}\}$$
is called the {\em covering radius} of $C$.

For $C$ a nonempty subset of $V\Gamma$ and for $0\leq i \leq\rho$, define
$$C_{i}=\{\hspace{1mm}x\in V\Gamma \ | \ d(x,C)=i\hspace{1mm}\}.$$
Then $\Pi(C)=\{C_{0}=C,C_{1},\ldots,C_{\rho}\}$ is the {\em distance
partition} of $V \Gamma$ with respect to code $C$.

A partition $\Pi=\{P_{0},P_{1},\ldots,P_{k}\}$ of $V\Gamma$
is called {\em equitable}
if, for all $i$ and $j$, the number of neighbors a vertex in
$P_{i}$ has in $P_{j}$ is independent of the choice of vertex in
$P_{i}.$
We say a code $C$ in $\Gamma$ is {\em completely regular} if this
distance partition $\Pi(C)$ is equitable\footnote{This definition of
a completely regular code is due to Neumaier \cite{G}. When $\Gamma$ is
distance-regular, it is equivalent to the
original definition, due to Delsarte \cite{D}, which we now mention. If
$\x$ is the characteristic vector of $C$, construct a
$|V\Gamma|\times(D+1)$ matrix with columns $A_i \x$ ($0\le i\le D$). Delsarte
declares $C$ to be completely regular if this {\em outer distribution matrix}
has only $\rho+1$ distinct rows.}. In this case the
following quantities are well-defined:
\begin{eqnarray}
\gamma_{i}= \left|\{y\in C_{i-1}\hspace{1mm}|\hspace{1mm}d(x,y)=1\} \right|,
\label{Egamma} \\
\alpha_{i}=\left|\{y\in C_{i}\hspace{1mm}|\hspace{1mm}d(x,y)=1\}\right|,
\hspace{4mm} \label{Ealpha} \\
\beta_{i}=\left|\{y\in C_{i+1}\hspace{1mm}|\hspace{1mm}d(x,y)=1\}\right| \hspace{3mm} \label{Ebeta}
\end{eqnarray}
where $x$ is chosen from $C_{i}$. The numbers $\gamma_{i},\alpha_{i},\beta_{i}$
are called the {\em intersection numbers} of code $C$. Observe that a graph
$\Gamma$ is distance-regular if and only if each vertex is a completely regular
code and these $|V\Gamma|$ codes all have the same intersection numbers.
An equitable partition
$\Pi=\{P_{1},\ldots,P_{m}\}$ of $V\Gamma$
is called a {\em completely regular partition} if all $P_i$ are completely regular codes
and any two of these have the same parameters.

If $\x$ is the characteristic vector of $C$ as a subset of $V\Gamma$,
then the {\em outer distribution module} of $C$ is defined as
$$\odm =\{M\x\hspace{1mm}|\hspace{1mm}M\in\bma\}.$$
Clearly, this is an $\bma$-invariant subspace of the
{\em standard module} $V:=\mathbb{C}^{V\Gamma}$.
Our next goal is to describe two nice bases for $\odm$.

For $0\le i \le \rho$, let $\x_i$ denote the characteristic vector of
$C_{i}$.

\begin{lemma}
\label{LodmA}
Let $\Gamma$ be a distance-regular graph and $C$ a completely regular code in
$\Gamma$. With notation as above,  we have
\begin{itemize}
\item[(a)] the vectors $\{ \x_0, \x_1, \ldots, \x_\rho \}$ form a basis for
the outer distribution module $\odm$ of $C$;
\item[(b)] relative to this basis, the matrix representing the action of
$A$ on $\odm$ is given by the tridiagonal matrix
$$U:=U(C)=\begin{pmatrix}
\alpha_{0} & \beta_{0}  &    &         &  \\
\gamma_{1} & \alpha_{1} & \beta_{1} &  &   \\
   & \gamma_{2} & \alpha_{2} &\beta_{2}  &   \\
   &            &  \ddots & \ddots  &  \ddots  \\
   &           & & \gamma_{\rho}&\alpha_{\rho}\end{pmatrix};$$
\item[(c)] $\dim \odm = \rho + 1$.
\end{itemize}
\end{lemma}

\begin{pf}
From Equations (\ref{Egamma}), (\ref{Ealpha}) and (\ref{Ebeta}) above, we have
\begin{equation}
 A \x_i =  \beta_{i-1} \x_{i-1} + \alpha_i \x_i + \gamma_{i+1} \x_{i+1} \label{EAxi}
\end{equation}
for $0\le i\le \rho$ where, for convenience, we set $\x_{-1}={\bf 0}$ and
$\x_{\rho+1}={\bf 0}$.
So a simple inductive argument shows that each $\x_i$ lies in the outer distribution
module of $C$. These vectors are trivially linear independent, so we need only
verify that they span $\odm$. By  (\ref{EAxi}), these vectors span an $A$-invariant
subspace of $V$ containing the characteristic vector $\x$ of $C$; since
$\odm$ is defined to be the smallest such subspace, the two spaces must coincide. \epf
\end{pf}

\begin{cor}
\label{CodmHadamardclosed}
Let $\Gamma$ be a distance-regular graph. For any completely regular code $C$
in $\Gamma$ with characteristic vector $\x$, the outer distribution
module $\odm$ of $C$ is closed under entrywise multiplication.
\end{cor}

\begin{pf}
Simply observe that the basis vectors $\x_i$ satisfy
$\x_i \circ \x_j = \delta_{i,j} \x_i$. \epf
\end{pf}

The tridiagonal matrix $U$ appearing in the lemma is called the {\em quotient matrix}
of $\Gamma$ with respect to $C$.

Now note that, for $0\le j\le D$, if the the vector $E_j \x$ is not the zero vector,
then it is an eigenvector for $A$ with eigenvalue $\theta_j$. This motivates
us to define
\begin{equation*}
S^*(C)=  \left\{ j \, | \, 1\le j\le D, \ E_j \x \neq {\bf 0} \right\}.
\end{equation*}

\begin{lemma}
\label{LodmE}
Let $\Gamma$ be a distance-regular graph and $C$ a completely regular code in
$\Gamma$. With notation as above,  we have
\begin{itemize}
\item[(a)] the nonzero vectors among the set
$\left\{ E_j \x \, | \, 0\le j\le D \right\}$ form a basis for
the outer distribution module $\odm$ of $C$;
\item[(b)] relative to this basis, the matrix representing the action of
$A$ on $\odm$ is a diagonal matrix with diagonal entries
$\left\{ \theta_j \, | \, j \in S^*(C) \cup \{0\} \right\}$;
\item[(c)]  $|S^*(C)| = \rho$.
\end{itemize}
\end{lemma}

\begin{pf}
Since $\bma$ is spanned both by
$\{ A_i\}_{i=0}^D$ and $\{ E_i\}_{i=0}^D$,
we see that $\odm$ is spanned by both
$\{ A_i\x \}_{i=0}^D$ and $\{ E_i \x\}_{i=0}^D$. Since the  nonzero vectors in
this latter set are linearly independent, they form a basis for $\odm$. From
Lemma \ref{LodmA}(c), we see that there must be exactly $\rho+1$ nonzero vectors in this
set, so $|S^*(C)| = \rho$. Finally, we have $A E_j \x = \theta_j E_j \x$ showing
that the matrix representing the action of $A$ on $\odm$ relative to this basis
is a diagonal matrix with diagonal entries as claimed. \epf
\end{pf}

\begin{cor}
\label{CeigU}
Let $\Gamma$ be a distance-regular graph and let  $C$ be a completely regular code
in $\Gamma$. With notation as above, the quotient matrix $U$ has $\rho+1$ distinct
eigenvalues, namely
$\left\{ \theta_j \, | \, j \in S^*(C) \cup \{0\} \right\}$.
\end{cor}

\begin{pf}
Suppose $S^*(C)=\{ i_1, \ldots , i_\rho\}$.
Since both $U$ and the diagonal matrix
$\diag  (k, \theta_{i_1}, \ldots, \theta_{i_\rho})$ represent
the same linear transformation, $A$, on the module $\odm$ with
respect to different bases, these two matrices must have the same
eigenvalues. \epf
\end{pf}

For $C$ a completely regular code in $\Gamma$,
we say that $\eta$ is an {\em eigenvalue of $C$} if
$\eta$ is an eigenvalue of the quotient matrix $U$ defined above. By
$\Ev(C)$, we denote the set of eigenvalues of $C$.  The above corollary
is often called ``Lloyd's Theorem'' in coding theory. The condition that each
eigenvalue of $C$ must be an eigenvalue of $\Gamma$ is a powerful condition
on the existence of completely regular codes, and perfect codes in particular\footnote{
A code $C$ in a distance-regular graph is {\em perfect} if $|C|=1$ or $\delta(C)
=2\rho(C)+1$. All perfect codes are completely rgeular.}.

Note that,
since $\gamma_i+\alpha_i+\beta_i=k$ for all $i$, $\theta_0=k$ belongs
to $\Ev(C)$. So
$$\Ev(C) = \{ k \} \cup \left\{ \theta_j \, | \, j \in S^*(C) \right\}.$$
Set ${\Ev}^{*}(C):=\Ev(C)-\{k\}$. For eigenvalue $\eta$ of $C$,
there is a unique right eigenvector
\begin{equation}
\label{Estdeigenvector}
u(\eta):=[u_{0}=1,u_{1},\ldots,u_{\rho}]^\top
\end{equation}
of $U$ associated to $\eta$; in analogy with the standard right eigenvectors
of graph $\Gamma$, we refer to this vector as the {\em standard (right)
eigenvector} of $C$ associated with $\eta$.
Note that this vector satisfies the following recurrence relation:
\begin{equation}
\label{Estdeigenvec}
\begin{split}u_{0}=1, \ u_{1}=
\frac{\eta-\alpha_{0}}{\beta_{0}},\hspace{25mm}\\
\gamma_{i}u_{i-1}+\alpha_{i}u_{i}+\beta_{i}u_{i+1}=\eta u_{i} \ \
(0\leq i \leq
\rho),\\ \text{where} \hspace{2mm} u_{-1}=u_{\rho+1}=0.\hspace{20mm}
\end{split}\end{equation}
For each standard right eigenvector of $C$, there is an eigenvector of $\Gamma$
in $\odm$
with the same eigenvalue which is unique up to scalar multiplication. For
eigenvalue $\theta_j$ of $C$, we refer to this {\em eigenvector belonging to $C$}
either as $E_j \x$ or as
\begin{equation} \label{Ebfu-u}
\mathbf{u} (\theta_j ) = \sum_{i=0}^\rho u_i \x_i
\end{equation}
where $u$ is defined above,
these two definitions differing only in their magnitude. Note that
$ \mathbf{u} (\theta_j ) \in \odm \cap V_j$.

\begin{lemma}
Assume that $\Gamma$ is $Q$-polynomial with $Q$-polynomial ordering
$\theta_{0}=k,\theta_{1},$ $\ldots,\theta_{D}$ of its eigenvalues.
Let $C$ be a completely regular code with
${\Ev}^{*}(C)=\{\theta_{i_{1}},\theta_{i_{2}},\ldots,\theta_{i_{\rho}}
\hspace{1mm}|\hspace{1mm}
i_{1}<i_{2}<\cdots<i_{\rho}\}$. Let
$\mathbf{u}(\theta_{i_{j}})$ be the eigenvector
with eigenvalue $\theta_{i_{j}}$ belonging to $C$. If
$\mathbf{u}(\theta_{i_{1}})$ has $\rho+1$ different entries, then
$i_{j}-i_{j-1}\leq i_{1}$ for all $j\in \{1,\ldots,\rho\}$.
\end{lemma}

\begin{pf}
By Lemma \ref{LodmA}(c), the outer distribution module
$\odm$ of $C$ has dimension $\rho+1$ and by Lemma \ref{LodmE}(a),
$ \left\{ E_j \mathbf{x} :  \theta_j \in \Ev(C) \right\}$ is a basis for it.
We now consider the entrywise product $\mathbf{u}^{(p)}$ of $p$ copies of the vector
$\mathbf{u} = E_{i_1} \x$.
Note that $\mathbf{u}^{(p)}\in \bma\mathbf{x}$ and that
$\Lambda:= \left\{ \mathbf{u}^{(p)}  :  0 \le p \le \rho \right\}$
is a linearly independent set of size $\rho+1$ by the Vandermonde property.
So $\Lambda$ spans $\odm$.
Suppose that $i_h - i_{h-1} \le i_1$ for $h<j$ but $i_j > i_{j-1}+i_1$.
Set
$$\mathcal{W}' = {\mbox {\rm span}} \left\{ E_0\mathbf{x}, E_{i_1}\mathbf{x},  \ldots,
 E_{i_{j-1}}\mathbf{x} \right\}.$$
As $\odm$ is closed under the Hadamard product, $ \mathbf{u} \circ \mathcal{W}'
\subseteq \odm $  and $ \mathbf{u} \circ \mathcal{W}' \subseteq V_0 + V_{i_{1}} + \cdots +
V_{i_{j-1}+i_1}$. Hence
$$ \mathbf{u} \circ \mathcal{W}' \subseteq \odm \cap (V_0 + V_{i_{1}} + \cdots + V_{i_{j-1}+i_1}).$$
But as $q_{i_1,h}^{l} =0$
for $h\le i_{j-1}$ and $l\ge i_j$, it follows  $ \mathbf{u} \circ \mathcal{W}'
\subseteq \mathcal{W}'$ and so $ \mathbf{u}^{(p)} \circ \mathcal{W}'
\subseteq \mathcal{W}'$ for $p\ge 1$ contradicting the fact that
$\Lambda$ spans $\odm$.\epf
\end{pf}

\begin{cor}
\label{Cgaps}
Let $\Gamma$ be a distance-regular graph and assume $\Gamma$
is $Q$-poly\-no\-mi\-al with respect to the natural ordering
$\theta_{0}=k>\theta_{1}>\cdots>\theta_{D}$ of its eigenvalues. Let $C$ be a
completely regular code in $\Gamma$ with
$S^*(C)=\{ i_1 , \ldots, i_\rho \} $
where $i_1 < \cdots < i_\rho$ and $\rho=\rho(C)$.
Then $i_{j}-i_{j-1}\leq i_{1}$ for all $j\in \{1,\ldots,\rho\}$.
\end{cor}
\begin{pf}
A standard argument involving Sturm sequences (see, e.g., {\cite[p.130]{B}} and
\cite[Lemma~8.5.2]{Godsil}) shows that, if $\theta_{i_1}$ is the second
largest eigenvalue of the tridiagonal matrix $U$, then the entries
of the standard right eigenvector of $C$ with respect to
$\theta_{i_1}$ are strictly decreasing. So the eigenvector $\mathbf{u}(\theta_{i_{1}})$
has $\rho+1$ distinct entries as required.\epf
\end{pf}

Our computational work suggests that Corollary \ref{Cgaps} is often a strong
feasibility condition for completely regular codes in the Hamming graphs.

Let $\Gamma$ be a distance-regular graph with diameter $D\geq 2$.
We say $\Gamma$ is an antipodal 2-cover whenever for all $x\in
V\Gamma$, there exists a unique vertex $y\in V\Gamma$ such that
$d(x,y)=D$. We denote this vertex by $\pi(x)$ and note that the
mapping $\pi:V\Gamma\longrightarrow V\Gamma$ is an
automorphism of $\Gamma$. It is known (cf. \cite[Prop.~4.2.3(ii)]{B})
that the subspace stabilized by this mapping is
$$ \left\{ {\bf v} \in V \ | \ {\bf v}_x = {\bf v}_{\pi(x)} \ \forall (x\in V\Gamma)
\right\}
= V_0 + V_2 + \cdots + V_{2\lfloor \frac{D}{2} \rfloor}  $$
and is therefore an $\bma$-submodule of the standard module.

\begin{lemma}
\label{L2cover}
Let $\Gamma$ be an antipodal 2-cover distance-regular graph and let
$\theta_{0}>\theta_{1}>\cdots>\theta_{D}$ be the distinct eigenvalues of
$\Gamma$. Let $C$ be a completely regular code with
$S^*(C) = \{ i_0 = 0 < i_1 < \ldots < i_\rho \}$
where $\rho=\rho(C)$. Let
$\pi$ be the automorphism defined above.
Then either
$$\pi(C)=C
\hspace{2mm}\text{and}\hspace{2mm} i_{j}\equiv 0
\hspace{3mm}(\text{mod}\hspace{1mm} 2)  \ \forall(j \in \{0,\ldots,\rho\})$$
or
$$\pi(C)=C_{\rho} \hspace{2mm} \text{and}\hspace{2mm} i_{j} \equiv j \hspace{3mm}(\text{mod}\hspace{1mm} 2)   \ \forall(j \in \{0,\ldots,\rho\}).$$
\end{lemma}

\begin{pf}
We know that $\odm$ is invariant under any $A_i$. So
$$ A_D \x = \tau_0 \x_0 + \cdots + \tau_\rho \x_\rho $$
for some scalars $\tau_0,\ldots,\tau_\rho$.
Let $x\in C$ and assume $\pi(x) \in C_{i}$ for some $i$.
Then $\tau_i \neq 0$ and so
for any vertex $y\in C_{i}$, $|\{z\in
C\hspace{1mm}|\hspace{1mm}d(y,z)=D\}|=1$.
This gives $C_{i}\subseteq\pi(C)$.  Since
$\rho(\pi(C))=\rho(C)$, the code $\pi(C)$ is either $C$ or $C_{\rho}$.

Let us first consider the case: $\pi(C)=C$.
In this case, the characteristic vector of $C$ belongs to the $A$-submodule
$ V_0 + V_2 + \cdots $ as outlined above, so for each $j$, $E_{i_j} \x$ belongs
to this submodule as well. Thus $i_j \equiv 0 \bmod{2}$ for all $0\le j\le \rho$.

In the other case, $\pi(C)=C_{\rho}$ and we use a Sturm sequence argument.
We know that $E_{i_j} \x$ is a scalar multiple of
$$ u_0 \x + u_1 \x_1 + \cdots + u_\rho \x_\rho $$
where $[u_0, u_1, \ldots, u_\rho]^\top$ is the standard eigenvector of $C$
associated with eigenvalue $\theta_{i_j}$.  But, by hypothesis, $\theta_{i_j}$
is the $j^{\rm th}$ largest eigenvalue of the tridiagonal quotient matrix $U$ defined
in the statement of Lemma \ref{LodmA}.  So by \cite[Lemma~8.5.2]{Godsil}, the
sequence $u_0,u_1,\ldots,u_\rho$ has $j$ sign changes. Since $u_0>0$, we
find $u_\rho$ is positive for $j$ even and negative for $j$ odd. But
it is well-known that  if $\vv$ is an eigenvector of an antipodal 2-cover
$\Gamma$, $\vv \in V_i$,
then $\vv_{\pi(x)} = \vv_x$ for  each $x \in V\Gamma$ when $i$ is even and
$\vv_{\pi(x)} = -\vv_x$ for  each $x \in V\Gamma$ when $i$ is odd. From this
we obtain our result. \epf
\end{pf}

\section{Q-polynomial Properties of a Code}
\label{Sec:Qpoly}
In this section, we will define $Q$-polynomial and Leonard completely regular codes and
establish a relation between them.

\begin{Def} \hfill
\label{DefQpoly}
Let $\Gamma$ be a distance-regular graph with diameter $D$ and \break
$\Ev(\Gamma)=\{\theta_{0}, \ldots,\theta_{i_{D}}\}.$
Let $C$ be a completely regular code with covering radius $\rho$ in $\Gamma.$
Then $C$ is called {\em $Q$-polynomial} if we have an ordering
${\Ev}(C)=\{\theta_{0}, \theta_{i_{1}},\ldots,\theta_{i_{\rho}}\}$ of the
eigenvalues of $C$ such that, for each $0\le p\le \rho$,
$\mathbf{u}^{(p)}:={\underbrace {\mathbf{u}\circ
\mathbf{u}\circ\cdots\circ \mathbf{u}}_{p \ \text{times}}}\in
\spn \{V_{i_{0}},\ldots,V_{i_{p}}\}$ where
$\mathbf{u} = E_{i_1} \mathbf{x} \in V_{i_{1}} $. In this case, we say $C$ is
{\em $Q$-polynomial with respect to} $\theta_{i_1}$.
\end{Def}

\begin{remark} Let $\Gamma$ be a distance-regular graph and
$x\in V\Gamma$. Then $C=\{x\}$ is completely regular and
$C$ is $Q$-polynomial with respect to the ordering $\theta_0,\theta_{i_1},
\ldots, \theta_{i_D}$ of $\Ev(C)$  if and only if
$\Gamma$ is $Q$-polynomial with respect to the ordering $E_0,E_{i_1},\ldots,
E_{i_D}$ of its primitive idempotents.
\end{remark}

Note that any completely regular code with covering radius at most 2 is
$Q$-polynomial. Also if we take for $C$ an antipodal pair in a
doubled Odd graph $\Gamma$ (see, for example \cite[Sec.~9.1D]{B}) then $C$ is $Q$-polynomial but $\Gamma$
is not $Q$-polynomial if its valency is at least 3.

Let $X$ be a finite abelian group. A {\em translation distance-regular graph}
on $X$ is a distance-regular graph $\Gamma$ with vertex set $X$ such that
if $x$ and $y$ are adjacent then $x+z$ and $y+z$ are adjacent for all $x,y,z\in X$.
A code $C\subseteq X$ is called {\em additive} for all $x,y\in C$, also $x-y \in C$;
i.e., $C$ is a subgroup of $X$.
If $C$ is an additive code in a translation distance-regular graph on $X$,
then we obtain the usual coset partition
$\Delta(C):=\{C+x \hspace{1mm}|\hspace{1mm} x\in X \}$
of $X$;  whenever $C$ is a completely regular code, it is easy to
see that $\Delta(C)$ is a completely regular partition. For any additive code
$C$ in a translation distance-regular graph $\Gamma$ on vertex set $X$, the
{\em coset graph} of $C$ in $\Gamma$ is the graph with vertex set $X/C$ and an edge
joining coset $C'$ to coset $C''$  whenever $\Gamma$ has an edge with one end in $C'$
and the other in $C''$. It follows from Theorem 11.1.6 in \cite{B} that this
coset graph is distance-regular whenever $C$ is an additive completely regular
code in a translation distance-regular graph.

\begin{prop}
\label{CQpolyquotient}
Let $X$ be a finite abelian group and let $\Gamma$ be a translation
distance-regular graph on $X$.  Let $C$ be an additive completely regular code in
$\Gamma$ and let $\Delta(C)$ be the partition of $X$ into cosets of $C$. Then
$C$ is $Q$-polynomial if and only if $\Gamma/\Delta(C)$ is a $Q$-polynomial
distance-regular graph.
\end{prop}

\begin{pf}
Let $C$ be an additive completely regular code in $\Gamma$ whose intersection numbers are $\gamma_i, \alpha_i$ and $\beta_{i}\hspace{1mm}(0\leq i\leq \rho)$. Then by \cite[p.352,353]{B}, eigenvalues of $\Gamma/\Delta(C)$ are $\frac{\eta_i-\alpha_0}{\gamma_1}$ for $\eta_i\in \Ev(C)$.
We see that $L(\Gamma/\Delta(C))=\frac{1}{\gamma_1}(U-\alpha_0I).$ Now the result follows easily.\epf
\end{pf}

\begin{Def}Let $\Gamma$ be a distance-regular graph. Let
$\eta$ be an eigenvalue of a completely regular code $C$ in
$\Gamma$ and let ${u}=[u_{0}=1,\ldots,u_{\rho}]^\top$ be the
standard eigenvector of $\eta$. Then the $\eta$ is called {\em
non-degenerate} if $u_{i-1}\neq u_{i}$ ($1\leq i\leq \rho$) and
$u_{i-1}\neq u_{i+1}$ ($1\leq i\leq \rho-1$).
\end{Def}

Note that the second largest eigenvalue of a completely regular
code is always non-degenerate. Likewise, if a code $C$ is $Q$-polynomial
with respect to the ordering $\{\eta_0,\eta_1,\ldots,\eta_\rho\}$ of its
eigenvalues, then $\eta_1$ is non-degenerate for $C$. This follows
from Definition \ref{DefQpoly}, which implies that the entrywise
powers of $\mathbf{u} =  E(\eta_1) \x$ are linearly independent and
Equation (\ref{Ebfu-u}) which then tells us that the $\rho+1$ entries of
the standard eigenvector for $\eta_1$ are all distinct.

\begin{prop}
\label{Plamtau}
Let $\Gamma$ be a distance-regular graph with valency $k$. Let $C$
be a completely regular code with covering radius $\rho$ and  $\Ev(C)=\{\eta_{i}\hspace{1mm}|\hspace{1mm}0\leq i\leq\rho\}$ in
$\Gamma$. Let
${u}(\eta_{i}):=[u_{0}=1,u_{1}(\eta_{i}),\ldots,u_{\rho}(\eta_{i})]^{T}$
be the standard eigenvector corresponding to eigenvalue
$\eta_{i}$ of $C$ $(0\leq i\leq\rho)$. Then there are (unique) $\lambda_{i},\tau_{i} \in \mathbb{R}$  such that $\sum_{i}\lambda_{i}=1,\sum_{i}\tau_{i}=1$ and the
following two hold:
\begin{equation}
\label{Eu-squared}
{u}^{(2)}(\eta_{1})=\sum_{i=0}^{\rho}\lambda_{i}{u}(\eta_{i})
\end{equation}
and
\begin{equation}
\label{Eu-cubed}
{u}^{(3)}(\eta_{1})=\sum_{i=0}^{\rho}\tau_{i} {u}(\eta_{i})
\end{equation}
In particular, if
$\eta_{1}$ is non-degenerate then the intersection numbers of $C$
are determined by the set of values
$$\{ \eta_0, \eta_1\} \cup \left\{\eta_{i} \hspace{1mm}|\hspace{1mm}\lambda_{i}\neq
0\hspace{2mm}\text{or}\hspace{2mm}\tau_{i}\neq 0 \right\}  \cup
\left\{ \lambda_0,\ldots, \lambda_\rho\right\} \cup \left\{ \tau_0,\ldots, \tau_\rho
\right\}. $$
\end{prop}

\begin{pf}
Let ${u}(\eta_{i})$ be the standard eigenvector of
$\eta_{i}$. The set
$\{{u}(\eta_{0}),\ldots,{u}(\eta_{\rho})\}$ forms a
basis of $\mathbb{R}^{\rho+1}$. Hence scalars $\lambda_i$ and $\tau_i$ each summing to one
and satisfying (\ref{Eu-squared}) and (\ref{Eu-cubed}) exist.

As
$\gamma_{j}u_{j-1}(\eta_{i})+\alpha_{j}u_{j}(\eta_{i})+\beta_{j}u_{j+1}(\eta_{i})=\eta_{i}u_{j}(\eta_{i}),$
(\ref{Eu-squared}) and (\ref{Eu-cubed}) can be rewritten as
$$\gamma_{j}u_{j-1}^{2}(\eta_{1})+\alpha_{j}u_{j}^{2}(\eta_{1})+\beta_{j}u_{j+1}^{2}(\eta_{1})=\sum_{i=0}^{\rho}\lambda_{i}\eta_{i}u_{j}(\eta_{i})$$
and
$$\gamma_{j}u_{j-1}^{3}(\eta_{1})+\alpha_{j}u_{j}^{3}(\eta_{1})+\beta_{j}u_{j+1}^{3}(\eta_{1})=\sum_{i=0}^{\rho}\tau_{i}\eta_{i}u_{j}(\eta_{i}).$$
Assume that we know the set $\{\eta_{i}
\hspace{1mm}|\hspace{1mm}\lambda_{i}\neq
0\hspace{2mm}\text{or}\hspace{2mm}\tau_{i}\neq
0\hspace{2mm}\text{or}\hspace{2mm}
i=0,1\}$ and all the $\lambda_i$ and $\tau_i$.
We use induction on $j$ to recover $\gamma_j,\alpha_j,\beta_j$ as well as
$u_{j+1}(\eta_i)$ for $1\le i\le \rho$. For $j=0$, the equations
$$\alpha_{0}+\beta_{0}=k,$$
$$\alpha_{0}+\beta_{0}u_{1}(\eta_{i})=\eta_{i}\hspace{5mm} \text{for}\hspace{2mm}0\leq i\leq\rho$$ and
$$\alpha_{0}+\beta_{0}u_{1}^{2}(\eta_{1})=\sum_{i=0}^{\rho}\lambda_{i}\eta_{i}.$$
easily allow us to obtain\footnote{Indeed, $\beta_0\neq 0$. If we
denote by $S$ the sum on the
right-hand side of the last equation, the simultaneous equations
$k+\beta_0( u_1(\eta_1)-1)=\eta_1$ and
$k+\beta_0( u_1(\eta_1)^2-1)=S$ allow us to solve for $u_1(\eta_1)+1$ and then for
$\beta_0$ so that all the remaining equations become linear.}
$\alpha_{0},\beta_{0},u_{1}(\eta_{i})$ for $0\leq i\leq\rho $.
Suppose that, for all $j\leq m$, the numbers
$\gamma_{j},\alpha_{j},\beta_{j},$ and $u_{j+1}(\eta_{i})$  ($0\leq i\leq\rho$) are
known. Now consider the case  $j=m+1$; we have four equations:
\begin{equation}
\label{Eind1}
\gamma_{m+1}+\alpha_{m+1}+\beta_{m+1}=k,
\end{equation}
\begin{equation}
\label{Eind2}
\gamma_{m+1}u_{m}(\eta_{1})+\alpha_{m+1}u_{m+1}(\eta_{1})+\beta_{m+1}u_{m+2}(\eta_{1})
=\eta_{1}u_{m+1}(\eta_{1}),
\end{equation}
\begin{equation}
\label{Eind3}
\gamma_{m+1}u_{m}^{2}(\eta_{1})+\alpha_{m+1}u_{m+1}^{2}(\eta_{1})+\beta_{m+1}u_{m+2}^{2}(\eta_{1})=\sum_{i=0}^{\rho}\lambda_{i}\eta_{i}u_{m+1}(\eta_{i})
\end{equation}
and
\begin{equation}
\label{Eind4}
\gamma_{m+1}u_{m}^{3}(\eta_{1})+\alpha_{m+1}u_{m+1}^{3}(\eta_{1})+\beta_{m+1}u_{m+2}^{3}(\eta_{1})=\sum_{i=0}^{\rho}\tau_{i}\eta_{i}u_{m+1}(\eta_{i}).
\end{equation}
As $\eta_{1}$ is non-degenerate, we obtain by Equations (\ref{Eind1})--(\ref{Eind4}):
\begin{eqnarray*}
u_{m+2}(\eta_1) \!\! &=& \!\! \frac{
R_\tau  - R_\lambda \left(u_{m+1}(\eta_1) + u_{m}(\eta_1) \right)
+ \eta_1 u_{m+1}^2(\eta_1)  u_{m}(\eta_1)
}{
R_\lambda + k u_{m+1}(\eta_1) u_{m}(\eta_1) - \eta_1 u_{m+1}(\eta_1)
\left(u_{m+1}(\eta_1) + u_{m}(\eta_1) \right)
}, \\[5mm]
\gamma_{m+1} &=& \! \frac{
R_\lambda + k u_{m+2}(\eta_1) u_{m+1}(\eta_1) - \eta_1 u_{m+1}(\eta_1)
\left(u_{m+2}(\eta_1) + u_{m+1}(\eta_1) \right)
}{
( u_{m}(\eta_1) - u_{m+2}(\eta_1) )(u_{m}(\eta_1) - u_{m+1}(\eta_1) ) }, \\[5mm]
\alpha_{m+1} &=& \frac{
R_\lambda + k u_{m+2}(\eta_1) u_{m}(\eta_1) - \eta_1 u_{m+1}(\eta_1)
\left(u_{m+2}(\eta_1) + u_{m}(\eta_1) \right)
}{
( u_{m+1}(\eta_1) - u_{m+2}(\eta_1) )(u_{m+1}(\eta_1) - u_{m}(\eta_1) ) }, \\[5mm]
\beta_{m+1} &=& \frac{
R_\lambda + k u_{m+1}(\eta_1) u_{m}(\eta_1) - \eta_1 u_{m+1}(\eta_1)
\left(u_{m+1}(\eta_1) + u_{m}(\eta_1) \right)
}{
( u_{m+2}(\eta_1) - u_{m+1}(\eta_1) )(u_{m+2}(\eta_1) - u_{m}(\eta_1) ) },
\end{eqnarray*}
where $R_\lambda$ and $R_\tau$ are shorthand for the expressions on the right-hand sides
of Equations (\ref{Eind3}) and (\ref{Eind4}), respectively; these quantities are
presumed known by the induction hypothesis.

But we also have, for $0\le i\le \rho$,
\begin{equation}
\label{Eind5}
\gamma_{m+1}u_{m}(\eta_{i})+\alpha_{m+1}u_{m+1}(\eta_{i})+\beta_{m+1}u_{m+2}(\eta_{i})=\eta_{i}u_{m+1}(\eta_{i})
\end{equation}
with (\ref{Eind1}) and (\ref{Eind2}) as special cases; from these,
we now obtain $u_{m+2}(\eta_{i})$ for $2\leq i\leq \rho$.  \epf
\end{pf}

\begin{lemma}
Let $\lambda_{j}$ and $\tau_{j}$ be the constants  defined in
Proposition~\ref{Plamtau} above. Suppose that
${\Ev}^{*}(C)=\{\theta_{i_{1}},\ldots,\theta_{i_{\rho}}\}$. If $\lambda_{j}\neq 0$, then
$q_{i_{1},i_{1}}^{i_{j}}\neq 0$ and if $\tau_{j}\neq 0$, then there
exists $i_{\ell}$ such that $q_{i_{1},{i_{1}}}^{i_{\ell}}\neq 0$ and
$q_{i_{\ell},i_{1}}^{i_{j}}\neq 0$.
\end{lemma}

\begin{pf}
Put $\mathbf{u}(\theta_{i_{j}}):= \sum_{h=0}^\rho u_h(\theta_{i_{j}})\mathbf{x}_h$.
Then
$\mathbf{u}^{(2)}(\theta_{i_{1}})=
\sum_{j=0}^{\rho}\lambda_{j}\mathbf{u}(\theta_{i_{j}})$,
$\mathbf{u}^{(3)}(\theta_{i_{1}})=
\sum_{j=0}^{\rho}   \tau_{j}\mathbf{u}(\theta_{i_{j}})$
and $\mathbf{u}(\theta_{i_{j}})\in V_{i_{j}}$. If
$\lambda_{j}\neq 0$ then as $\mathbf{u}^{(2)}(\theta_{i_{1}})$
is not orthogonal to $V_{i_{j}}$, by
Theorem~\ref{Tcgs}, $q_{i_{1},i_{1}}^{i_{j}}\neq 0$.
Since
$\mathbf{u}^{(3)}(\theta_{i_{1}})=\sum_{\ell=0}^{\rho}\lambda_{\ell}
\mathbf{u}(\theta_{i_{\ell}})\circ\mathbf{u}(\theta_{i_{1}})$,
if $\tau_{j}\neq 0$ then there exists $\ell$ such that
$\lambda_{\ell}\neq 0$ and
$\mathbf{u}(\theta_{i_{\ell}})\circ\mathbf{u}(\theta_{i_{1}})$ is not
orthogonal to $V_{i_{j}}$, by  Theorem~\ref{Tcgs}, there exists
$i_{\ell}$ such that $q_{i_{1},i_{1}}^{i_{\ell}}\neq 0$ and
$q_{i_{\ell},i_{1}}^{i_{j}}\neq 0.$ \epf
\end{pf}

Let $\Gamma$ be a distance-regular graph with adjacency matrix $A$
and let $C\subseteq V\Gamma$ be a completely regular code with
covering radius $\rho$, ${\Ev}^{*}(C)=\{\theta_{i_{1}},\ldots,$ $\theta_{i_{\rho}}\}$ and distance partition
$\{C_{0}, C_{1},\ldots,C_{\rho}\}$. For $0\leq i\leq\rho$, let
$\mathbf{x}_{i}$ denote the characteristic vector of subconstituent $C_{i}$.
Let $\mathcal{B}^{*}
:=\{\mathbf{x}_{i}\hspace{1mm}|\hspace{1mm}i=0,\ldots,\rho\}$ and
$\mathcal{B}:=\{E_{i_{j}}\mathbf{x}_{0}\hspace{1mm}|\hspace{1mm}j=0,\ldots,\rho\}$.
Then both $\mathcal{B}^{*}$ and $\mathcal{B}$ are
bases for the outer distribution module $\odm$ of $C$.
Now consider first the linear
transformation $\mathcal{A}$ on $\odm$ which is defined by
$\mathcal{A}(\mathbf{y})=A\mathbf{y}$ for $\mathbf{y}\in
\odm$. For any nontrivial eigenvalue $\theta$ of $C$,
define the linear transformation $\mathcal{A}^{*}(\theta)$ on
$\odm$ by
$\mathcal{A}^{*}(\theta)(\mathbf{y})=(E(\theta)\mathbf{x}_{0})\circ\mathbf{y}$
for $\mathbf{y} \in \odm$. Since
$A\mathbf{x}_{i}=\beta_{i-1}\mathbf{x}_{i-1}+\alpha_{i}\mathbf{x}_{i}+\gamma_{i+1}\mathbf{x}_{i+1}$,
the matrix representing $\mathcal{A}$ with respect to the basis
$\mathcal{B}^{*}$ is irreducible tridiagonal (i.e. each entry on
the subdiagonal and each entry on the superdiagonal are nonzero)
and the matrix representing $\mathcal{A}$ with respect to the
basis $\mathcal{B}$ is diagonal. We can easily check that the
matrix representing $\mathcal{A}^{*}(\theta)$ with respect to the
basis $\mathcal{B}^{*}$ is diagonal as
$(E(\theta)\mathbf{x}_{0})\circ\mathbf{x}_{i}=(E(\theta)\mathbf{x}_{0})_{y}\mathbf{x}_{i}$
where $y\in C_{i}$. We now define a Leonard completely regular
code.

\begin{Def}
With above notation, a completely regular code $C$ is called {\em
Leonard} if there exists a nontrivial eigenvalue $\theta$ of $C$
such that the matrix representing
$\mathcal{A}^{*}=\mathcal{A}^{*}(\theta)$ with respect to
$\mathcal{B}$ is irreducible tridiagonal. When this happens for a particular
eigenvalue $\theta$, we will say that $C$ is {\em Leonard with respect to $\theta$.}
 Note that, following Terwilliger {\em \cite[p.150]{J}}, the pair
$\mathcal{A},\mathcal{A}^{*}$ is a Leonard pair on $\odm$ for a Leonard completely regular code.
\end{Def}

\begin{prop}
\label{Pleonard}
Let $\Gamma$ be a distance-regular graph. Then any Leonard
completely regular code of $\Gamma$ is a $Q$-polynomial completely
regular code.
\end{prop}

\begin{pf}
Let $C$ be a completely regular code with covering radius
$\rho$ and characteristic vector $\x$.
Suppose $C$ is Leonard with respect to  the nontrivial eigenvalue $\theta$ of $C$.
Since the matrix representing $\mathcal{A}^{*}(\theta)$
is irreducible tridiagonal with respect
to some ordering of the basis $\mathcal{B}$,
we may index   ${\Ev}^{*}(C)=\{\theta_{i_{1}},\ldots,\theta_{i_{\rho}}\}$
so that,  for $0\leq j\leq \rho$ we have
$E_{i_{1}}\mathbf{x} \circ E_{i_{j}}\mathbf{x} =\epsilon_{j}
E_{i_{j-1}}\mathbf{x} +\varphi_{j}
E_{i_{j}}\mathbf{x} +\psi_{j} E_{i_{j+1}}\mathbf{x}$ for
some scalars $\epsilon_{j},\varphi_{j},\psi_{j}$ ($\epsilon_j$ and $\psi_j$ being nonzero) where
$E_{i_{-1}}\mathbf{x} = E_{i_{\rho+1}}\mathbf{x} =0$.
The result follows. \epf
\end{pf}

\begin{Def}
We say a Leonard code is of type {\em Krawtchouk} if the corresponding 
Leonard pair is of type {\em Krawtchouk} as defined in Terwilliger 
{\em \cite{L}}. In a similar fashion, we define Leonard codes of type
{\em Hahn}, {\em dual Hahn}, {\em Racah} and so on.\\
Sometimes we also say that a Leonard code is of {\em class} 
(\textrm{I}), (\textrm{IA}), (\textrm{IB}), (\textrm{II}), (\textrm{IIA}), 
(\textrm{IIB}), (\textrm{IIC}), (\textrm{IID}) and (\textrm{III}) if 
the corresponding Leonard pair is of class (\textrm{I}), (\textrm{IA}),
(\textrm{IB}), (\textrm{II}), (\textrm{IIA}), (\textrm{IIB}), 
(\textrm{IIC}), (\textrm{IID}) and (\textrm{III}), respectively, where 
we use the notation of Bannai and Ito {\em \cite{A}}.
\end{Def}

It is a natural problem to choose one of these families and to classify
all Leonard codes of that type. It is interesting to note that
a  Leonard code of a given type may appear within a
classical distance-regular graph of some other type.  For example,
the $n$-cube is obviously a $Q$-polynomial distance-regular graph of
Krawtchouk type, and it contains the binary repetition code, which is
not of Krawtchouk type. Below, in Example \ref{Erifa}, we describe additive
binary completely regular codes found by Rifa and Zinoviev which are of
dual Hahn type.

Let $\theta$ be a eigenvalue of $C$ and
$\mathcal{A}^{*}:=\mathcal{A}^{*}(\theta)$. For $0\leq i\leq\rho$,
as
$\mathcal{A}^{*}\mathbf{x}_{i}=(E(\theta)\mathbf{x})_{y}\mathbf{x}_{i}$
where $y\in C_{i}$, the vector $\mathbf{x}_{i}$ is an eigenvector
for $\mathcal{A}^{*}$. Let $F_{j}^{*}$ and $F_{j}$ denote the
primitive idempotent corresponding to $\mathbf{x}_{j}$ and
$E_{i_{j}}\mathbf{x}$, respectively. In \cite[Lemma 5.7]{K},
Terwilliger shows that if at least three of the following four
conditions hold then $\mathcal{A},\mathcal{A}^{*}$ is a Leonard
pair.
\begin{equation}
\label{Eu-tri1}
F_{h}^{*}\mathcal{A}F_{j}^{*} \ \begin{cases} = 0 \hspace{8mm}\text{if} \hspace{2mm} h-j>1\\
\neq0 \hspace{4mm}\text{if} \hspace{2mm} h-j=1
\end{cases}
\hspace{3mm}(0\leq h,j\leq \rho),
\end{equation}
\begin{equation}
\label{Eu-tri2}
F_{h}^{*}\mathcal{A}F_{j}^{*} \ \begin{cases} = 0 \hspace{8mm}\text{if} \hspace{2mm} j-h>1\\
\neq0 \hspace{4mm}\text{if} \hspace{2mm} j-h=1
\end{cases}
\hspace{3mm}(0\leq h,j\leq \rho),
\end{equation}
\begin{equation}
F_{h}\mathcal{A}^{*}F_{j} \ \begin{cases} = 0 \hspace{8mm}\text{if} \hspace{2mm} h-j>1\\
\neq 0 \hspace{4mm}\text{if} \hspace{2mm} h-j=1
\end{cases}
\hspace{3mm}(0\leq h,j\leq \rho),
\end{equation}
\begin{equation}
F_{h}\mathcal{A}^{*}F_{j} \ \begin{cases} = 0 \hspace{8mm}\text{if} \hspace{2mm} j-h>1\\
\neq 0 \hspace{4mm}\text{if} \hspace{2mm} j-h=1
\end{cases}
\hspace{3mm}(0\leq h,j\leq \rho).
\end{equation}
Note that Equations (\ref{Eu-tri1}) and (\ref{Eu-tri2}) together imply that
the matrix representing
$\mathcal{A}$ with respect to $\mathcal{B}^{*}$ is irreducible tridiagonal.

\begin{theorem}
\label{TQiffLeonard}
\hfill Let $\Gamma$ be a distance-regular graph with diameter $D$ and \\
${\Ev}(\Gamma)=\{\theta_{0},\ldots,\theta_{D}\}$. Let $C$ be a completely
regular code in $\Gamma$. Then $C$ is Leonard if and only if C is
$Q$-polynomial.
\end{theorem}

\begin{pf}
The `only if' part is done by Proposition \ref{Pleonard}, so we only need to
show the `if' part. Let $C$ be a completely regular
code with covering radius $\rho$ which is $Q$-polynomial with respect to 
eigenvalue $\theta$.  Definition \ref{DefQpoly} then gives us a natural ordering
${\Ev}^{*}(C)=\{\theta_{i_{1}},\ldots,\theta_{i_{\rho}}\}$ where $\theta_{i_1}=\theta$.
Let $\mathcal{A}^{*}:=\mathcal{A}^{*}(\theta)$ and we now
consider the products $F_{h} \mathcal{A}^{*} F_{j}$
for $0\leq h,j\leq \rho$. As
$\mathcal{A}^{*}E_{i_{j}} \x  = E_{i_{1}} \x \circ
E_{i_j} \x $ and as $C$ is $Q$-polynomial, there exists a
polynomial $p_{j+1}$ of degree exactly $j+1$ such that
$\mathcal{A}^{*}E_{i_{j}} \x  = p_{j+1}(E_{i_{1}} \x).$
Since $\mathcal{B}$ is a basis for $\odm$ and 
$\mathcal{A}^{*}E_{i_{j}} \x  \in \odm$, we can write
$\mathcal{A}^{*}E_{i_{j}} \x =\sum_{l=0}^{\rho}\xi_{l}E_{i_{l}} \x $
where $\xi_{l}\in\mathbb{R} \ (0\leq l\leq\rho)$ satisfy the following condition:
$$\xi_{l} \ \begin{cases} = 0 \hspace{8mm}\text{if} \hspace{2mm} l>j+1\\
\neq 0 \hspace{4mm}\text{if}\hspace{2mm} l=j+1 \end{cases}.$$
Observe
$F_{h}E_{i_{l}} \x =\delta_{h,l}E_{i_{l}} \x $
for $0\leq h,l\leq \rho$. By this, we find
$$F_{h}\mathcal{A}^{*}E_{i_{j}} \x  \ \begin{cases} = 0 \hspace{8mm}\text{if} \hspace{2mm} h>j+1\\
\neq 0 \hspace{4mm}\text{if}\hspace{2mm} h=j+1\end{cases}.
$$
So 
$$F_{h}\mathcal{A}^{*}F_{j} \ 
        \begin{cases} = 0 \hspace{8mm}\text{if} \hspace{2mm} h-j>1\\
                   \neq 0 \hspace{4mm}\text{if}\hspace{2mm} h-j=1\end{cases}.
$$  
and the result follows. \epf
\end{pf}

In \cite{L}, Terwilliger gave a parametrization of any Leonard pair. It
follows that, for any Leonard pair, there are at most seven free parameters.
(Allowing for equivalence under affine transformations, this may be reduced to
five.) We now show that the Leonard pair associated to a $Q$-polynomial
completely regular code in a known distance-regular graph has all its
parameters determined by just six free parameters.

\begin{cor}
\label{CQpolyfp}
Let $\Gamma$ be a distance-regular graph of valency $k$ and diameter
$D$.  Let $C$ be a completely regular code in $\Gamma$ which is $Q$-polynomial
with respect to the ordering $\eta_0,\eta_1,\ldots, \eta_\rho$ of $\Ev(C)$.
Then the 
intersection numbers $\alpha_i,\beta_i,\gamma_i$ ($0\le i\le \rho$) are
completely determined (as is the covering radius $\rho$, from $\beta_\rho=0$)
by the eigenvalues $\eta_1$ and $\eta_2$
of $C$ together with the parameters $\lambda_0$, $\lambda_1$,
$\tau_1$ and $\tau_2$ as defined in Proposition \ref{Plamtau}.
\end{cor}

\begin{pf}
Since  $C$ is $Q$-polynomial, we have
\begin{equation}
\label{Eu-squared-Q}
{u}^{(2)}(\eta_{1})=   \lambda_{0}{u}(\eta_{0}) +  \lambda_{1}{u}(\eta_{1}) +
\lambda_{2}{u}(\eta_{2})
\end{equation}
and
\begin{equation}
\label{Eu-cubed-Q}
{u}^{(3)}(\eta_{1})= \tau_0 {u}(\eta_0) + \tau_1 {u}(\eta_1) +
\tau_2 {u}(\eta_2) + \tau_3 {u}(\eta_3) .
\end{equation}
Looking at the zero entry on both sides of each equation, we find
$\lambda_0+\lambda_1+\lambda_2=1$ and $\tau_0+\tau_1+\tau_2+\tau_3=1$.  Now $C$ is
Leonard by Theorem \ref{TQiffLeonard}, so there exist scalars
$\sigma_1,\sigma_2,\sigma_3$ for which
\begin{equation}
\label{Eu-u2-Q}
{u}(\eta_{1}) \circ u(\eta_2) = \sigma_1 {u}(\eta_1) + \sigma_2 {u}(\eta_2) + \sigma_3 {u}(\eta_3) .
\end{equation}
Moreover, we have   $\sigma_1+\sigma_2+\sigma_3 = 1$.  Next,
we may use this and Equation (\ref{Eu-squared-Q}) to obtain
an alternative expression for $u^{(3)}(\eta_1)$:
$$
{u}^{(3)}(\eta_{1})= \lambda_0 \lambda_1 {u}(\eta_0) + \left(\lambda_0+ \lambda_1^2 + \lambda_2 \sigma_1\right) {u}(\eta_1)
+ \lambda_2 \left( \lambda_1+\sigma_2 \right) {u}(\eta_2) + \lambda_2 \sigma_3 {u}(\eta_3) .
$$
Comparing coefficients against those in Equation (\ref{Eu-cubed-Q}), we find
\begin{eqnarray*}
\lambda_0 \lambda_1  &=& \tau_0  \\
\lambda_0 + \lambda_1^2 +  \lambda_2 \sigma_1 &=& \tau_1 \\
\lambda_2 \left( \lambda_1+\sigma_2 \right)  &=& \tau_2 \\
\lambda_2 \sigma_3 &=& \tau_3
\end{eqnarray*}
so that $\lambda_2,\tau_0,\tau_3$ are determined by knowledge of
$\lambda_0$, $\lambda_1$, $\tau_1$ and $\tau_2$. Now all we need are
the eigenvalues needed in Proposition \ref{Plamtau}. But we know $\eta_0=k$,
the valency of $\Gamma$, we are given $\eta_1$ and $\eta_2$ by hypothesis
and we may then solve for $\eta_3$ by looking at the $i=1$  entry on both sides
of (\ref{Eu-cubed-Q}):
$$ \tau_0 + \tau_1 \frac{ \eta_1 - \alpha_0 }{ k - \alpha_0 } +
\tau_2 \frac{ \eta_2 - \alpha_0 }{ k - \alpha_0 }  +
\tau_3 \frac{ \eta_3 - \alpha_0 }{ k - \alpha_0 } =
\left( \frac{ \eta_1 - \alpha_0 }{ k - \alpha_0 } \right)^3$$
where we have used  the evaluation (\ref{Estdeigenvec})
$u_1(\theta)= (\theta-\alpha_0)/(k-\alpha_0)$. Now the result follows from
Proposition \ref{Plamtau}. \epf
\end{pf}

\begin{conj}
Every completely regular code in a $Q$-polynomial distance-regular graph
with sufficiently large covering radius is a Leonard completely regular code.
\end{conj}

We finish this section with a description of an interesting family of
codes in the $n$-cubes.

\begin{Ex}
\label{Erifa}
In any $\binom{m}{2}$-cube for integer $m\geq 3$, there exist
Leonard completely regular codes  which are not of Krawtchouk type.
Following \cite{H}, for natural numbers $m\geq 3$ and $2\le l < m$,
define $E_{l}^{m}$ as the set of all binary vectors of length $m$ and weight $l$.
Denote by $H^{(m,l)}$ the binary matrix of size
$m\times\binom{m}{l}$, whose columns are exactly all vectors from
$E_{l}^{m}$. Rifa and Zinoviev consider the binary linear code
$C^{(m,l)}$ whose parity check matrix is the matrix $H^{(m,l)}$; they
show that the code $C^{(m,2)}$ is completely regular and its coset
graph is the halved $m$-cube. As the halved $m$-cube is $Q$-polynomial,
it follows that $C^{(m,2)}$ is Leonard, but it is of dual Hahn, not Krawtchouk, type.
\end{Ex}

\section{Harmonic completely regular codes}
\label{Sec:harmonic}

In a companion paper \cite{klm+},  we explore a well-structured class of
Leonard completely regular codes in the Hamming graphs. These {\em arithmetic}
completely regular codes are defined as those whose eigenvalues are in
arithmetic progression: $\Ev(C)= \{ k, k-t,k-2t,\ldots\}$. These codes have a
rich structure and are intimately tied to Hamming quotients of Hamming graphs.
In \cite{klm+}, we study products of completely regular codes
and completely classify the possible quotients of a Hamming graph that can
arise from the coset partition of a linear arithmetic completely regular code.
For families of distance-regular graphs other than the Hamming graphs, we
need to look at a slightly weaker definition to probe the same sort of
rich structure.

We next introduce the class of {\em harmonic} completely regular codes and we
will see that this class lies strictly between the arithmetic completely regular codes
and the Leonard completely regular codes.

\begin{Def}
Let $\Gamma$ be a $Q$-polynomial distance-regular graph with
respect to the ordering $\theta_{0},\theta_{1},\ldots,\theta_{D}$
of its eigenvalues and $C$ be a completely regular code of
$\Gamma$. We call the code $C$ {\em harmonic} if
$\Ev(C)=\{\theta_{ti}\hspace{1mm}|\hspace{1mm}i=0,\ldots,\rho\}$ for some
positive integer $t$.
\end{Def}

Let $\Gamma$ be a $Q$-polynomial with respect to the ordering
$\{\theta_{0},\theta_{1},\ldots,\theta_{D}\}$ of its eigenvalues and
let $C\subseteq V\Gamma$ be a code.
Then {\em strength} of $C$, $t(C)$ is defined as the $\min\{i\geq 1
\hspace{1mm}|\hspace{1mm}\theta_{i}\in {\Ev}^*(C)\}-1.$
\begin{Ex}
The following are examples of  harmonic completely regular
codes:\\
(1) the repetition code in a hypercube;\\
(2) cartesian products of a completely regular code of a Hamming
graph $C\times\cdots \times C$ where $C$ is covering radius $1$;\\
(3) in the Grassmann Graph $J_{q}(n,t)$, whose vertices are all $t$-dimensional
subspaces of a some $n$-dimensional vector space $V$ over $GF(q)$, we find
the following two families:
\begin{itemize}
\item $C$ consists of all $t$-dimensional subspaces of a given $(n-s)$-dimensional
subspace of $V$, where $0<s<n-t$;
\item  $C$ consists of all $t$-dimensional subspaces of $V$ containing a fixed
$s$-dimensional subspace $U$ of $V$, where $0<s \le t < n$.
\end{itemize}
(We note that the Johnson graph $J(n,t)$ contains examples analogous to these.) \hfill\break
(4) any completely regular code of strength $0$ in a
$Q$-polynomial distance-regular graph.
\end{Ex}

\begin{lemma}
Let $\Gamma$ be a $Q$-polynomial distance-regular graph with
respect to the ordering $\theta_{0},\theta_{1},\ldots,\theta_{D}$
of its eigenvalues. Then any harmonic completely regular code is
a Leonard completely regular code.
\end{lemma}
\begin{pf} \hfill
Since $\Gamma$ is $Q$-polynomial, there exist numbers
$\omega_{h,j}$ such that \\
$E_{t}\mathbf{x}_{0}\circ
E_{jt}\mathbf{x}_{0}=\sum_{h=0}^{\rho}\omega_{h,j}E_{ht}\mathbf{x}_{0}$
and the following holds:
$$\omega_{h,j} \ \begin{cases} = 0 \hspace{8mm}\text{if} \hspace{2mm}
|ht-jt|>t\\
\neq 0 \hspace{4mm}\text{if}\hspace{2mm}|ht-jt|\leq
t\end{cases}.$$ So,
$$\omega_{h,j} \   \begin{cases} = 0 \hspace{8mm}\text{if} \hspace{2mm}
|h-j|>1\\
\neq 0 \hspace{4mm}\text{if}\hspace{2mm}|h-j|\leq 1\end{cases}.$$
Hence the matrix representing $\mathcal{A}^{*}(\theta_{t})$ is
irreducible tridiagonal with respect to $\mathcal{B}$. \epf
\end{pf}

Finally, we remark that the codes given in Example \ref{Erifa} are Leonard
but not harmonic.

\section*{Acknowledgments}

Part of this work was completed while the third author was visiting Pohang
Institute of Science and Technology (POSTECH). WJM wishes to thank the
Department of Mathematics at POSTECH for their hospitality and Com$^2$MaC
for financial support.
JHK and LWS are partially supported by
the Basic Science Research Program through the National Research
Foundation of Korea (NRF) funded by the Ministry of Education,
Science and Technology (grant number 2009-0089826). JHK was also partially
supported by a grant of the Korea Research Foundation
funded by the Korean Government (MOEHRD) under grant number
KRF-2007-412-J02302.
WJM wishes to thank the US National Security Agency for financial support
under grant number H98230-07-1-0025.

The authors wish to thank
Paul Terwilliger for helpful discussions regarding some of the material in this
paper.

\end{document}